\documentclass[11pt,twoside,a4paper]{article}
\usepackage{amssymb}
\usepackage{amsmath}
\usepackage[dvips]{graphicx}
\usepackage{amscd}
\title{Algebraic Shifting and Basic Constructions on Simplicial Complexes}
\author{Eran Nevo \footnote{Institute of Mathematics, The Hebrew
University, Givat Ram, Jerusalem 91904, Israel,
e-mail address: eranevo@math.huji.ac.il}}
 \newtheorem{thm}{Theorem}[section]
 \newtheorem{cor}[thm]{Corollary}
 
\newtheorem{de}[thm]{Definition}
 \newtheorem{prop}[thm]{Proposition}
\newtheorem{prob}[thm]{Problem}
\newtheorem{conj}[thm]{Conjecture}

\begin{document}
\maketitle
\begin{abstract}
We try to understand the behavior of algebraic shifting with
respect to some basic constructions on simplicial complexes, such
as union, coning, and (more generally) join. In particular, for
the disjoint union of simplicial complexes we prove
$\Delta(K\dot{\cup}L)=\Delta(\Delta(K)\dot{\cup}\Delta(L))$
(conjectured by Kalai \cite{skira}), and for the join we give an
example of simplicial complexes $K$ and $L$ for which
$\Delta(K*L)\neq\Delta(\Delta(K)*\Delta(L))$ (disproving a
conjecture by Kalai \cite{skira}), where $\Delta$ denotes the
(exterior) algebraic shifting operator. We develop a 'homological'
point of view on algebraic shifting which is used throughout this
work.
\end{abstract}

\section{Introduction}\label{sec0}
 Algebraic shifting is an operator which associates with each simplicial
complex another simplicial complex which is combinatorially
simpler (is shifted: to be defined shortly) and preserves some
properties of the original complex. It was introduced by Kalai
\cite{55}.
\\

  In this work we try to understand the behavior of algebraic
shifting with respect to some basic constructions on simplicial
complexes, such as union, cone, and (more generally) join. In some
cases we get a 'nice' behavior:
We prove that the disjoint union of simplicial complexes satisfies
$$\Delta(K\dot{\cup}L)=\Delta(\Delta(K)\dot{\cup}\Delta(L))$$
(conjectured by Kalai \cite{skira}). Moreover, we give an explicit
combinatorial description of $\Delta(K\dot{\cup}L)$ in terms of
$\Delta(K)$ and $\Delta(L)$.
%characterization of the $3$-tuples
%$$(\Delta(K),\Delta(L),\Delta(K\dot{\cup}L))$$
These results follow from the following theorem about shifting a
(not necessarily disjoint) union: Define $init_{j}(S)$ to be the
set of lexicographically least $j$ elements in $S$, for
$S\subseteq \mathbb{N}$, where $\mathbb{N}$ is the set of positive
integers endowed with the usual order. For every $i>0$ define
 $I_{S}^{i} = I_{S}^{i}(n) = \{T: T\subseteq [n], |T|=|S|+i, init_{|S|}(T)=S\}$. Denote by $K_{i}$ the
$i$-th skeleton of a simplicial complex $K$.

\begin{thm} \label{gen}
Let $K$ and $L$ be two simplicial complexes, and let $d$ be the
dimension of $K\cap L$. For every
   subset $A$ of the set of vertices $[n]=(K\cup L)_{0}$, the
   following additive formula holds:
   \begin{eqnarray}\label{eq.genKuL}
|I_{A}^{d+2}\cap \Delta(K\cup L)|=|I_{A}^{d+2}\cap
\Delta(K)|+|I_{A}^{d+2}\cap \Delta(L)|.
   \end{eqnarray}
\end{thm}
In the case where $K\cap L$ is a simplex (with all of its
subsets), the following stronger assertion holds:
\begin{thm}\label{clique-sum}
Let $K$ and $L$ be two simplicial complexes, where $K\cap
L=<\sigma>$ is the simplicial complex consisting of the set
$\sigma$ and all of its subsets. For every $i>0$ and every
   subset $S$ of the set of vertices $[n]=(K\cup L)_{0}$, the
   following additive formula holds:
$$|I_{S}^{i}\cap \Delta(K\cup_{\sigma}L)|=|I_{S}^{i}\cap\Delta(K)|+|I_{S}^{i}\cap\Delta(L)|-|I_{S}^{i}\cap<\sigma>|.$$

\end{thm}
This theorem gives an explicit combinatorial description of
$\Delta(K\cup L)$ in terms of $\Delta(K)$, $\Delta(L)$ and
$dim(\sigma)$. In particular, any gluing of $K$ and $L$ along a
$d$-simplex results in the same shifted complex $\Delta(K\cup L)$,
depending only on $\Delta(K)$, $\Delta(L)$ and $d$.

Another instance of a 'nice' behavior is
$$\Delta\circ Cone =
Cone \circ\Delta$$ (Kalai \cite{skira}). We prove a generalized
version of this property for near-cones (defined in [1]). In the
case of join, we do not get such a good behavior: We give an
example of simplicial complexes $K$ and $L$ for which
$$\Delta(K*L)\neq\Delta(\Delta(K)*\Delta(L))$$ where $*$ denotes the
join operator, disproving a conjecture by Kalai \cite{skira}.
%Although this property fails even for suspension (join with two
%points), when restricting to the faces of maximal dimension which
%represent the homology, it holds.
However, a weaker assertion holds:
%For every $i$, The number of
%faces of maximal dimension in $K*L$ which do not intersect $[i]$
%equals the product of the numbers of such faces in $K$ and in $L$.
\begin{thm}\label{K*L,max}
Let $K$ and $L$ be two simplicial complexes with $|(K*L)_{0}|=n$.
Then for every $i\in [n]$:
\begin{eqnarray}\label{join}
|\{S\in \Delta(K*L): [i]\cap S=\emptyset,
|S|=dim(K*L)+1\}|=\nonumber\\
|\{S\in \Delta(K): [i]\cap S=\emptyset, |S|=dim(K)+1\}| \times
\nonumber\\
|\{S\in \Delta(L): [i]\cap S=\emptyset, |S|=dim(L)+1\}|.
\end{eqnarray}
\end{thm}
The case $i=1$ was known - it follows from the K\"{u}nneth theorem
with field coefficients \cite{Mun} and from the combinatorial
interpretation of homology using shifting (Bj\"{o}rner and Kalai
\cite{BK}). We give an example that Theorem \ref{K*L,max} does not
hold for symmetric shifting.
\\
For a survey on algebraic shifting, the reader can consult
\cite{skira}.
\\
Outline: Section \ref{sec0+} introduces the needed algebraic
background. Section \ref{sec1} develops a 'homological' point of
view on the algebraic shifting operator which will be used in the
successive sections. Section \ref{sec2} shifts unions of
simplicial complexes (proving Theorems \ref{gen} and
\ref{clique-sum}), Section \ref{sec4} shifts near-cones, Section
\ref{sec k*l} shifts joins of simplicial complexes (proving
Theorem \ref{K*L,max}).

\section{Algebraic Background}\label{sec0+}
   In this section we set some needed preliminaries.
We follow the definitions and notation of \cite{skira}. Let $K$ be
a simplicial complex on a vertex set $[n]$. The $i$-th skeleton of
$K$ is $K_{i}=\{S\in K: |S|=i+1\}$. For each $1\leq k \leq n$ let
$<_{L}$ be the lexicographic order on $(^{[n]}_{\ k})$, i.e.
$S<_{L}T \Leftrightarrow min\{a:a\in S\triangle T\}\in S$, and let
$\triangleleft_{P}$ be the partial order defined as follows: Let
$S=\{s_{1}<\dots<s_{k}\}, T=\{t_{1}<\dots<t_{k}\}$,
$S\triangleleft_{P} T$ iff $s_{i}\leq t_{i}$ for every $1\leq
i\leq k$ ($min$ and $\leq$ are taken with respect to the usual
order on $\mathbb{N}$).

Let $V$ be an $n$-dimensional vector space over $\mathbb{R}$ with
basis $\{e_{1},\dots,e_{n}\}$. Let $\bigwedge V$ be the graded
exterior algebra over $V$. Denote
$e_{S}=e_{s_{1}}\wedge\dots\wedge e_{s_{j}}$ where $S=
\{s_{1}<\dots<s_{j}\}$. Then $\{e_{S}: S\in (_{\ j}^{[n]})\}$ is a
basis for $\bigwedge^{j} V$. Note that as $K$ is a simplicial
complex, the ideal $(e_{S}:S\notin K)$ of $\bigwedge V$ and the
vector subspace $span\{e_{S}:S\notin K\}$ of $\bigwedge V$ consist
of the same set of elements in $\bigwedge V$. Define the exterior
algebra of $K$ by
$$\bigwedge (K)=(\bigwedge V)/(e_{S}:S\notin K).$$ Let
$\{f_{1},\dots,f_{n}\}$ be a basis of $V$, generic over
$\mathbb{Q}$ with respect to $\{e_{1},\dots,e_{n}\}$, which means
that the entries of the corresponding transition matrix $A$
($e_{i}A=f_{i}$ for all $i$) are algebraically independent over
$\mathbb{Q}$. Let $\tilde{f}_{S}$ be the image of $f_{S}\in
\bigwedge V$ in $\bigwedge(K)$. Define
$$\Delta(K)=\Delta_{A}(K)=\{S: \tilde{f}_{S}\notin span\{\tilde{f}_{S'}:S'<_{L}S\}\}$$
to be the shifted complex (introduced by Kalai \cite{55}).
 The construction is canonical, i.e. it is independent of the choice of the generic matrix
 $A$, and for a permutation $\pi:[n]\rightarrow [n]$ the induced simplicial complex
  $\pi(K)$ satisfies $\Delta(\pi(K))=\Delta(K)$.
 It results in a shifted simplicial complex, having the
 same face vector and Betti vector as $K$'s \cite{BK}
(for a simplicial complex $L$, its face vector $f(L)=(f_{i}(L))$
is defined by $f_{i}(L)=|L_{i}|$, and its Betti (homology) vector
$\beta(L)=(\beta_{i}(L))$ is defined by $\beta_{i}(L)=dim
\tilde{H}_{i}(L,\mathbb{R})$, where $\tilde{H}_{i}()$ stands for
the reduced $i$-th homology). The key ingredient in Bj\"{o}rner
and Kalai's proof that algebraic shifting preserves Betti numbers,
is the following combinatorial way of reading them:
$$\beta_{i}(L)=|\{S\in \Delta(L)_{i}: S\cup 1 \notin \Delta(L)\}|.$$

Fixing the basis $\{e_{1},\dots,e_{n}\}$ of $V$ induces the basis
$\{e_{S}: S\subseteq [n]\}$ of $\bigwedge V$, which in turn
induces the dual basis $\{e_{T}^{*}: T\subseteq [n]\}$ of
$(\bigwedge V)^{*}$ by defining $e_{T}^{*}(e_{S})=\delta_{T,S}$.
($(\bigwedge V)^{*}$ stands for the space of $\mathbb{R}$-linear
functionals.) For $f,g\in \bigwedge V$ $<f,g>$ will denote
$f^{*}(g)$. Define the so called left interior product of $g$ on
$f$ \cite{56}, where $g,f \in \wedge V$, denoted $g\lfloor f$, by
the requirement that for all $h\in\bigwedge V$
$$<h,g\lfloor f>=<h\wedge g,f>.$$
($g\lfloor\cdot$ is the adjoint operator of $\cdot\wedge g$ w.r.t.
the inner product $<\cdot,\cdot>$ on $\bigwedge V$.) Thus,
$g\lfloor f$ is a bilinear function, satisfying
\begin{equation}\label{floorE}
e_{T}\lfloor e_{S}=\{^{\pm e_{S\backslash T}\  if\  T\subseteq S}
_{0 \ otherwise}
\end{equation}
where the sign equals $(-1)^{a}$, where $a=|\{(s,t)\in S\times T:
s\notin T, t<s\}|$. This implies in particular that for a monomial
$g$ (i.e. $g$ is a wedge product of elements of degree 1) $g
\lfloor$ is a boundary operation on $\bigwedge V$, and in
particular on $span\{e_{S}: S\in K\}$ \cite{56}.

Let us denote for short $\bigwedge^{j}K=span\{e_{S}: S\in
K_{j-1}\}$ and $\bigwedge K=\bigwedge(K)=span\{e_{S}: S\in K\}$.
This should cause no confusion with the definition of the exterior
algebra of $K$, as we shall never use this exterior (quotient)
algebra structure in the sequel, but only the graded vector space
structure $\bigwedge K$ just defined.
 We denote:
$$Ker_{j}f_{R}\lfloor(K) = Ker_{j}f_{R}\lfloor=Ker (f_{R}\lfloor: \bigwedge^{j+1}K \rightarrow \bigwedge^{j+1-|R|}K).$$
Note that the definition of $\bigwedge K$ makes sense more
generally when $K_{0}\subseteq [n]$ (and not merely when
$K_{0}=[n]$), and still $f_{R}\lfloor$ operates on the subspace
$\bigwedge K$ of $\bigwedge V$ for every $R\subseteq [n]$. (Recall
that $f_{i}=\alpha_{i1}e_{1}+...+\alpha_{in}e_{n}$ where
$A=(\alpha_{ij})_{1\leq i,j\leq n}$ is a generic matrix.) Define
$f^{0}_{i}=\sum_{j\in K_{0}}\alpha_{ij}e_{j}$, and
$f^{0}_{S}=f^{0}_{s_{1}}\wedge\dots\wedge f^{0}_{s_{j}}$ where $S=
\{s_{1}<\dots<s_{j}\}$. By equation (\ref{floorE}), the following
equality of operators on $\bigwedge K$ holds:
\begin{equation} \label{floorF}
\forall S\subseteq [n] \ \ f_{S}\lfloor=f^{0}_{S}\lfloor.
\end{equation}
We now turn to find relations between the kernels defined above
and algebraic shifting.

\section{Shifting and Kernels of Boundary Operations}\label{sec1}
In this section we give some equivalent descriptions of the
algebraic shifting operator, using the kernels defined in Section
\ref{sec0+}. This approach will be used throughout this work.
\\
 The following generalizes a result for graphs \cite{56} (the
proof is similar):
\begin{prop} \label{prop.1}
Let $R\subseteq [n]$, $|R|<j+1$. Then: $$Ker_{j}f_{R}\lfloor =
\bigcap_{i:[n]\ni i\notin R}Ker_{j}f_{R\cup{i}}\lfloor.$$
\end{prop}
%\begin{proof}
$Proof$: Recall that $h \lfloor(g\lfloor f)=(h\wedge g)\lfloor f$.
Thus, if $f_{R}\lfloor m=0$ then
$$f_{R\cup i}\lfloor m=\pm f_{i}\lfloor (f_{R}\lfloor
m)=f_{i}\lfloor 0=0.$$ Now suppose $m\in span\{e_{S}: S\in
K,|S|=j+1\}\setminus Ker_{j}f_{R}\lfloor$. The set $\{f_{Q}^{*}:
Q\subseteq[n],|Q|=j+1-|R|\}$ forms a basis of $(\bigwedge
^{j+1-|R|}V)^{*}$,
 so there
is some $f_{R'}$, $R'\subseteq[n], |R'|=j+1-|R|$, (note that
$R'\neq \emptyset$) such that
$$<f_{R'}\wedge f_{R},m>=<f_{R'},f_{R}\lfloor m>\neq0.$$
We get that for $i_{0}\in R'$: $i_{0}\notin R$ and $<f_{R'\setminus i_{0}},f_{R\cup
i_{0}}\lfloor m>\neq 0 $. Thus $m\notin Ker_{j}f_{R\cup i_{0}}\lfloor$ which completes
the proof.$\square$
%\end{proof}
\\
In the next proposition we determine the shifting of a simplicial
complex by looking at the intersection of kernels of boundary
operations (actually only at their dimensions): Let $S$ be a
subset of $[n]$ of size $s$. For $R\subseteq [n]$,$|R|=s$, we look
at $f_{R}\lfloor:\bigwedge^{s}(K)\rightarrow
\bigwedge^{s-|R|}(K)=\mathbb{R}$.
%Proposition \ref{prop.2} is the
%"dual" of the previous definition of shifting in the sense that
%$f\lfloor$ and $f\wedge$ are dual operations.
\begin{prop} \label{prop.2}
  Let $K_{0}, S\subseteq[n], |K_{0}|=k, |S|=s$. The following
  quantities are equal:
\begin{eqnarray}\label{KerShift}
dim \bigcap_{R<_{L}S, |R|=s, R\subseteq [n]}Ker_{s-1}f_{R}\lfloor,\label{KS1}\\
dim \bigcap_{R<_{L}S, |R|=s, R\subseteq [k]}Ker_{s-1}f^{0}_{R}\lfloor,\label{KS2}\\
|\{T\in \Delta(K): |T|=s, S\leq_{L} T\}|.\label{KS3}
\end{eqnarray}
In particular, $S\in \Delta(K)$ iff $$dim\bigcap_{R<_{L}S, |R|=s,
R\subseteq [n]}Ker_{s-1}f_{R}\lfloor > dim\bigcap_{R\leq_{L}S,
|R|=s, R\subseteq [n]}Ker_{s-1}f_{R}\lfloor$$ (equivalently, $S\in
\Delta(K)$\ iff\ \  $\bigcap_{R<_{L}S, |R|=s, R\subseteq
[n]}Ker_{s-1}f_{R}\lfloor \nsubseteq Ker_{s-1}f_{S}\lfloor$).
\end{prop}
%\begin{proof}
$Proof$:
%Here $f_{R}\lfloor$ is into $\mathbb{R}$.
First we show that (\ref{KS1}) equals (\ref{KS2}). For every
$T\subseteq [n]$, $T<_{L}S$, decompose $T=T_{1}\cup T_{2}$,
 where $T_{1}\subseteq [k], T_{2}\cap [k]=\emptyset$. For each
 $T_{3}$ satisfying $T_{3}\subseteq [k], T_{3}\supseteq T_{1},
 |T_{3}|=|T|$,
 we have $T_{3}\leq_{L}T$. Each $f^{0}_{t}$, where $t\in T_{2}$, is a linear
 combination of the $f^{0}_{i}$'s, $1\leq i \leq k$, so $f^{0}_{T}$ is a
 linear combination of such $f^{0}_{T_{3}}$'s. Thus, for every
 $j\geq s-1$,
 $$\bigcap_{R\leq_{L}T,R\subseteq
 [k]}Ker_{j}f^{0}_{R}\lfloor \subseteq Ker_{j}f^{0}_{T}\lfloor,$$ and
 hence $$\bigcap_{R<_{L}S,R\subseteq
 [k]}Ker_{j}f^{0}_{R}\lfloor = \bigcap_{R<_{L}S,R\subseteq
 [n]}Ker_{j}f^{0}_{R}\lfloor.$$ Combining with (\ref{floorF}) the desired equality follows.

Next we show that (\ref{KS3}) equals (\ref{KS2}). Let $m\in
\bigwedge^{s}(K)$ and $R\subseteq[n], |R|=s$. Let us express $m$
and $f_{R}$ in the basis $\{e_{S}: S\subseteq [n]\}$:
$$m=\sum_{T\in K,|T|=s} \gamma_{T}e_{T}$$
$$f_{R}=\sum_{S'\subseteq[n],|S'|=s}A_{RS'}e_{S'}$$
where $A_{RS'}$ is the minor of $A$ (transition matrix) with
respect to the rows $R$ and columns $S'$, and where $\gamma_{T}$
is a scalar in $\mathbb{R}$.

By bilinearity we get $$f^{0}_{R}\lfloor m =f_{R}\lfloor m
=\sum_{T\in K,|T|=s}\gamma_{T}A_{RT}.$$ Thus
% $dim\bigcap_{R<_{L}S,|R|=S}Ker_{S-1}f_{R}\lfloor$
(\ref{KS2}) equals the dimension of the solution space of the
system $B_{S}x=0$, where $B_{S}$ is the matrix $(A_{RT})$, where
$R<_{L}S$, $R\subseteq [k]$, $|R|=s$ and $T\in K,|T|=s$. But,
since the row indices of $B_{S}$ are an initial set with respect
to the lexicographic order, the intersection of $\Delta(K)$ with
this set of indices determines a basis of the row space of
$B_{S}$. Thus, $rank(B_{S})=|\{R\in \Delta(K): |R|=s,R<_{L}S\}|$.
But $K$ and $\Delta(K)$ have the same $f$-vector, so we get:
$$dim \bigcap_{R<_{L}S, R\subseteq [k], |R|=s}Ker_{s-1}f^{0}_{R}\lfloor =
f_{s-1}(K) - rank(B_{S})=$$ $$= |\{T\in \Delta(K):|T|=s, S\leq_{L}
T\}|$$ as desired.$\square$
%\end{proof}
\\
  Recall that for each $j>0$ and $S\subseteq [n]$, $|S|\geq j$ we define $init_{j}(S)$ to be the set of lexicographically least $j$
  elements in $S$, and for every $i>0$
 $I_{S}^{i} = I_{S}^{i}(n) = \{T: T\subseteq [n], |T|=|S|+i, init_{|S|}(T)=S\}$. Let
  $S_{(i)}^{(m)} = S_{(i)}^{(m)}(n)=min_{<_{L}}I_{S}^{i}(n)$ and
 $S_{(i)}^{(M)} = S_{(i)}^{(M)}(n)=max_{<_{L}}I_{S}^{i}(n)$.
In the sequel, all the sets of numbers we consider are subsets of
$[n]$. In order to simplify notation, we will often omit noting
that. We get the following information about the partition of the
faces in the shifted complex into 'intervals':
\begin{prop} \label{prop.*}
Let $K_{0}\subseteq [n]$, $S\subseteq [n]$, $i>0$. Then:
$$|I_{S}^{i}\cap \Delta(K)| =  dim \bigcap_{R<_{L}S}Ker_{|S|+i-1}f_{R}\lfloor(K) -
dim \bigcap_{R\leq_{L}S}Ker_{|S|+i-1}f_{R}\lfloor(K).$$
\end{prop}
$Proof$:
   By Proposition \ref{prop.1},$$dim \bigcap_{R<_{L}S}Ker_{|S|+i-1}f_{R}\lfloor=
   dim \bigcap_{R<_{L}S}\bigcap_{j\notin R, j\in [n]}Ker_{|S|+i-1}f_{R\cup j}\lfloor =...=$$
$$dim \bigcap_{R<_{L}S}\ \bigcap_{T:T\cap R=\emptyset, |T|=i}Ker_{|S|+i-1}f_{R\cup
T}\lfloor =$$
   $$dim \bigcap_{R<_{L}S_{(i)}^{(m)}}Ker_{|S_{(i)}^{(m)}|-1}f_{R}\lfloor$$
(to see that the last equation is true, one needs to check that
$\{R\cup T: T\cap R=\emptyset, |T|=i, R<_{L}S\} = \{Q:
Q<_{L}S_{(i)}^{(m)}\}$). By Proposition \ref{prop.2},$$dim
\bigcap_{R<_{L}S_{(i)}^{(m)}}Ker_{|S_{(i)}^{(m)}|-1}f_{R}\lfloor =
|\{Q\in\Delta(K):|Q|=|S|+i,S_{(i)}^{(m)}\leq_{L}Q\}|.$$ Similarly,
$$dim \bigcap_{R\leq_{L}S}Ker_{|S|+i-1}f_{R}\lfloor(K)=|\{F\in\Delta(K):|F|=|S|+i,S_{(i)}^{(M)}<_{L}F\}|$$
 (here one checks
that $\{R\cup T: T\cap R=\emptyset, |T|=i, R\leq_{L}S\} = \{F:
F\leq_{L}S_{(i)}^{(M)}\}$). Thus, the proof of the proposition is
completed. $\square$
\\
Note that on $I_{S}^{1}$ the lexicographic order and the partial
order $\triangleleft_{P}$ coincide, since all sets in $I_{S}^{1}$
have the same $|S|$ least elements. As $\Delta(K)$ is shifted,
$I_{S}^{1}\cap\Delta(K)$ is an initial set of $I_{S}^{1}$ with
respect to $<_{L}$. Denote for short
$$D(S)=D_{K}(S)=|I^{1}_{init_{|S|-1}(S)}(n)\cap \Delta(K)|.$$
$D_{K}(S)$ is indeed independent of the particular $n$ we choose,
as long as $K_{0}\subseteq [n]$. We observe that:
\begin{prop} \label{prop.3}
Let $K_{0}$ and $S=\{ s_{1}<\dots<s_{j}<s_{j+1}\}$ be subsets of
$[n]$. Then: $S\in \Delta(K)\Leftrightarrow s_{j+1}-s_{j}\leq
D(S).$ $\square$
\end{prop}

Another easy preparatory lemma is the following:
\begin{prop} \label{prop.5}
Let $K_{0}, S\subseteq [n]$. Then: $D_{\Delta K}(S)=D_{K}(S).$
\end{prop}
%\begin{proof}
$Proof$: It follows from the fact that $\Delta^{2}=\Delta$ (Kalai
\cite{D^2=D}, or later on here in Corollary \ref{prop.13
delta^2}). $\square$
%\end{proof}

\section{Shifting Union of Simplicial Complexes}\label{sec2}
Let us consider a general union first:
\begin{prob}\label{U problem}(\cite{skira}, Problem 13)
Given two simplicial complexes $K$ and $L$, find all possible
connections between $\Delta(K \cup L)$, $\Delta(K)$, $\Delta(L)$
and $\Delta(K \cap L)$.
\end{prob}
 We look on $\bigwedge(K\cup L)$, $\bigwedge(K\cap
L)$, $\bigwedge(K)$ and
  $\bigwedge(L)$ as subspaces of $\bigwedge(V)$ where
  $V=span\{e_{1},...,e_{n}\}$ and $[n]=(K\cup L)_{0}$. As before,
  the $f_{i}$'s are generic linear combinations of the $e_{j}$'s
  where $j\in [n]$. Let $S\subseteq [n]$, $|S|=s$ and $1\leq j$.
First we find a connection between boundary operations on the
spaces associated with $K$, $L$, $K\cap L$ and $K\cup L$ via the
following commutative diagram of exact sequences:
\begin{equation} \label{diagram}
\begin{CD}
0@> >>\bigwedge^{j+s}(K\cap L)@>i>>\bigwedge^{j+s}K \bigoplus \bigwedge^{j+s}L@>p>>\bigwedge^{j+s}(K\cup L)@> >>0\\
@VV V @VVfV @VVgV @VVhV @VV V \\
0@> >>\oplus\bigwedge^{j}(K\cap L)@>\oplus i>>\oplus\bigwedge^{j}K \bigoplus \oplus\bigwedge^{j}L@>\oplus p>>\oplus\bigwedge^{j}(K\cup L)@> >>0\\
\end{CD}
\end{equation}
where all sums $\oplus$ in the bottom sequence are taken over
$\{A: A<_{L}S, |A|=s\}$ and $i(m)=(m,-m)$, $p((a,b))=a+b$, $\oplus
i(m)=(m,-m)$, $\oplus p((a,b))=a+b$,
$f=\oplus_{A<_{L}S}f_{A}\lfloor(K\cap L)$,
$g=(\oplus_{A<_{L}S}f_{A}\lfloor(K),\oplus_{A<_{L}S}f_{A}\lfloor(L))$
and $h=\oplus_{A<_{L}S}f_{A}\lfloor(K\cup L)$.

By the snake lemma, (\ref{diagram}) gives rise to the following
exact sequence:
\begin{equation} \label{snake}
\begin{CD}
0@> >>ker f@> >>ker g@> >>ker h@>\delta>> \\
coker f@> >>coker g@> >>coker h@> >>0\\
\end{CD}
\end{equation}
where $\delta$ is the connecting homomorphism. Let
($\ref{diagram}'$) be the diagram obtained from (\ref{diagram}) by
replacing $A<_{L}S$ with $A\leq_{L}S$ everywhere, and renaming the
maps by adding a superscript to each of them. Let ($\ref{snake}'$)
be the sequence derived from ($\ref{diagram}'$) by applying to it
the snake lemma. If $\delta=0$ in (\ref{snake}), and also the
connecting homomorphism $\delta'=0$ in ($\ref{snake}'$), then by
Proposition \ref{prop.*} the following additive formula holds:
\begin{equation} \label{additive-formula}
|I_{S}^{j}\cap \Delta(K\cup L)|=|I_{S}^{j}\cap
\Delta(K)|+|I_{S}^{j}\cap \Delta(L)|-|I_{S}^{j}\cap \Delta(K\cap
L)|.
\end{equation}
\subsection{A proof of Theorem \ref{gen}}\label{cup1}
$Proof\  of\  Theorem \ \ref{gen}$: Put $j=d+2$ in (\ref{diagram})
and in ($\ref{diagram}'$). Thus, the range and domain of $f$ in
(\ref{diagram}) and of $f'$ in ($\ref{diagram}'$) are zero, hence
$ker f=coker f=0$ and $ker f'=coker f'=0$, and Theorem \ref{gen}
follows. $\square$
\\

\textbf{Remark}:  It would be interesting to understand what extra
information about $\Delta(K\cup L)$ we can derive by using more of
the structure of $\Delta(K\cap L)$, and not merely its dimension.
In particular, it would be interesting to find combinatorial
conditions that imply the vanishing  of $\delta$ in (\ref{snake}).
The proof of Theorem \ref{clique-sum} in Subsection
\ref{KsimplexL} provides a step in this direction.
 The Mayer-Vietoris long exact sequence (see \cite{Mun} p.186) gives some information of this type,
  by the interpretation of the Betti vector using the shifted complex [1], mentioned in Section \ref{sec0+}.

\subsection{How to shift a disjoint union?}\label{disj}
As a corollary to Theorem \ref{gen} we get the following
combinatorial formula for shifting the disjoint union of
simplicial complexes:
\begin{thm} \label{prop.6}
%(our main result in this section).
  Let $(K\dot{\cup}L)_{0}=[n], [n]\supseteq S=\{ s_{1}<\dots<s_{j}<s_{j+1}\}$. Then: $$S\in\Delta(K\dot{\cup}L)
  \Leftrightarrow s_{j+1}-s_{j}\leq |I^{1}_{init_{|S|-1}(S)}\cap\Delta(K)|+|I^{1}_{init_{|S|-1}(S)}\cap\Delta(L)|.$$
\end{thm}
%\begin{proof}

 $Proof$: Put $d=-1$ and $A=init_{|S|-1}(S)$ in Theorem \ref{gen},
 and by Proposition \ref{prop.3} we are done.
 $\square$
%\end{proof}
\\

  As a corollary, we get the following nice equation, proposed by
Kalai \cite{skira}:

\begin{cor} \label{prop7}
  $\Delta(K\dot{\cup}L)=\Delta(\Delta(K)\dot{\cup}\Delta(L)).$
\end{cor}
%\begin{proof}
$Proof$: $S\in\Delta(K\cup L)$ iff (by Theorem \ref{prop.6}) $s_{j+1}-s_{j}\leq
D_{K}(S)+D_{L}(S)$ iff (by Proposition \ref{prop.5}) $s_{j+1}-s_{j}\leq
D_{\Delta(K)}(S)+D_{\Delta(L)}(S)$ iff (by Theorem \ref{prop.6}) $S\in
\Delta(\Delta(K)\dot{\cup}\Delta(L))$. $\square$
%\end{proof}
\\
\textbf{Remarks}: (1) Above a high enough dimension (to be
specified) all faces of the shifting of a union are determined by
the shifting of its components: Let $st(K\cap L)=\{\sigma\in K\cup
L: \sigma\cap(K\cap L)_{0}\neq \emptyset \}$. Then $\Delta(K)$ and
$\Delta(L)$ determine all faces of $\Delta(K\cup L)$ of dimension
$>dim(st(K\cap L))$, by applying Theorem \ref{prop.6} to the
subcomplex of $K\cup L$ spanned by the vertices $(K\cup
L)_{0}-(K\cap L)_{0}$.

    (2) Let $X$ be a $(k+l)\times (k+l)$
generic block matrix, with an upper block of size $k\times k$ and
a lower block of size $l\times l$. Although we defined the
shifting operator $\Delta=\Delta_{A}$ with respect to a generic
matrix $A$, the definition makes sense for any nonsingular matrix
(but in that case the resulting complex may not be shifted). Let
$K_{0}=[k]$ and $L_{0}=[k+1,k+l]$. Corollary \ref{prop7} can be
formulated as
$$\Delta\circ\Delta_{X}(K\dot{\cup}L)=\Delta(K\dot{\cup}L)$$
because $\Delta_{X}(K\dot{\cup}L)=\Delta(K)\dot{\cup}\Delta(L)$
(on the right hand side of the equation the vertices of the two
shifted complexes are considered as two disjoint sets). However,
there are simplicial complexes $C$ on $k+l$ vertices , for which
$\Delta\circ\Delta_{X}(C)\neq\Delta(C)$. For example, let $k=l=3$
and take the graph $G$ of the octahedron
$\{\{1\},\{4\}\}*\{\{2\},\{5\}\}*\{\{3\},\{6\}\}$. Then
$\Delta\circ\Delta_{X}(G)\ni \{4,5\} \notin \Delta (G)$.

    (3) By induction, we get from Corollary \ref{prop7} that:
$$\Delta(\dot{\cup}_{1\leq i\leq n}K^{i})=\Delta(\dot{\cup}_{1\leq i\leq n}\Delta(K^{i}))$$
for any positive integer $n$ and disjoint simplicial complexes
$K^{i}$.

  (4) Theorem \ref{prop.6} gives a very simple (linear time in $t=2^{n}$) algorithm for
computing $\Delta(K\dot{\cup}L)$, given $\Delta(K)$ and
$\Delta(L)$, $n=|(K\cup L)_{0}|$.

  (5) For symmetric algebraic shifting, introduced by Kalai
  \cite{KalaiSymm}, the results about shifting the disjoint union,
  namely the symmetric analogous assertions to Theorem \ref{prop.6} and Corollary
  \ref{prop7}, remain true.
   Techniques similar to those developed in
  Sections \ref{sec1} and in the proof of Theorem \ref{gen} provide a proof of
  the symmetric analogue of Theorem \ref{gen} in the case where
  $d=-1$. We omit the details.
  For a general union, we do not know
  whether the symmetric analogue of Theorem \ref{gen} holds or not.

\subsection{A recursive formula for shifting a disjoint union}
    We now turn to prove a recursive formula for computing
$\Delta(K \dot{\cup} L)$, given $\Delta(K)$ and $\Delta(L)$,
conjectured by Kalai \cite{skira}.

We introduce the  $\sqcup$ operator (as in \cite{skira}) defined
on shifted simplicial complexes: Let $K,L$ be shifted simplicial
complexes. Define $K\sqcup L$ to be the simplicial complex (yet to
be shown) satisfying $(K\sqcup L)_{0}=K_{0}\dot{\cup} L_{0}=[k+l]$
($K_{0}=[k], L_{0}=[l]$, but regarding $K_{0}, L_{0}$ as disjoint
sets) and recursively satisfying
$$(K\sqcup L)_{j}= 1*(lk(1,K)\sqcup lk(1,L))_{j-1}\cup(ast(1,K)\sqcup
ast(1,L))_{j}$$ where the link and anti-star of $S\in K$ are
defined by: $$lk(S,K)=\{T\in K: T\cap S=\emptyset,T\cup S\in
K\},$$
$$ast(S,K)=\{T\in K: T\cap S=\emptyset\}.$$
Note that $lk(S,K)\subseteq ast(S,K)$. The $1*()$ operator on sets
of sets is defined (unusually) as: $1*S=\{1\dot{\cup} s: s\in S\}$
for $S$ such that for all $s\in S$ $1\notin s$.  Kalai
conjectured:

\begin{cor} \label{prop.8}
Let $K,L$ be two shifted simplicial complex. Then:
$$\Delta(K\dot{\cup}L)=K\sqcup L.$$
\end{cor}
%\begin{proof}
$Proof$: By definition $(K\sqcup L)_{0}=\Delta(K\dot{\cup}L)_{0}$.
We proceed by induction on $n=|(K\dot{\cup}L)_{0}|$. By the
induction hypothesis, for every $j$
$$(K\sqcup L)_{j}=[1*\Delta(lk(1,K)\dot{\cup}
lk(1,L))_{j-1}]\cup[\Delta(ast(1,K)\dot{\cup} ast(1,L))_{j}].$$
Denoting
$$Ast=ast(1,K)\dot{\cup}ast(1,L),$$
$$Lk=lk(1,K)\dot{\cup}lk(1,L),$$
we have $$K\sqcup L=1*\Delta(Lk)\cup\Delta(Ast).$$
 First we show that $K\sqcup L$ is a simplicial complex:
$Lk\subseteq Ast$ are two simplicial complexes, therefore
$\Delta(Lk)\subseteq\Delta(Ast)$, and we get that $K\sqcup L$ is
also a simplicial complex.
 Our second step is to show that $K\sqcup L$ is shifted.
As $1*\Delta(Lk)$ (not a simplicial complex) and $\Delta(Ast)$ are
shifted, we only have to show that
$$\partial(\Delta(Ast))\subseteq\Delta(Lk).$$ (For a set $A\subseteq
(^{[n]}_{\ k})$, $\partial(A)=\{b\in(^{\ [n]}_{k-1}): \exists a\in
A\  such that\  b\subseteq a\}$. For a simplicial complex $K$,
$\partial(K)=\bigcup_{i}
\partial(K_{i})$.) A basic property of algebraic shifting \cite{skira} is that for every simplicial complex
$C$,
$$\partial(\Delta(C))\subseteq\Delta(\partial(C)),$$
so we get $$\partial(\Delta(Ast))\subseteq\Delta(\partial(Ast)).$$
As $K$ and $L$ are shifted, $\partial(ast(1,K))\subseteq lk(1,K)$,
and similarly for $L$. Hence $\partial(Ast)\subseteq Lk$, and
therefore $\Delta(\partial(Ast))\subseteq\Delta(Lk)$. Now that we
know that $K\sqcup L$ is a shifted simplicial complex, by
Proposition \ref{prop.3} it is enough to show that for each
$S\subseteq[n]$, we have
$$D_{K\sqcup L}(S)=D_{\Delta(K\dot{\cup}L)}(S),$$ as $(K\sqcup L)\cap
I^{1}_{init_{|S|-1}(S)}$ is an initial set of
$I^{1}_{init_{|S|-1}}(S)$.
\\*
case 1: $1\notin S$: For each $T\in I^{1}_{init_{|S|-1}}(S)$: $1\notin T$,
therefore
$$D_{K\sqcup L}(S)=D_{\Delta(Ast)}(S)=D_{ast(1,K)}(S)+D_{ast(1,L)}(S)$$
$$=D_{K}(S)+D_{L}(S)=D_{\Delta(K\dot{\cup}L)}(S).$$
The second and last equations are by Theorem \ref{prop.6} and
Proposition \ref{prop.5}.
\\*
case 2: $1\in S$: For each $T\in I^{1}_{init_{|S|-1}}(S)$: $1\in T$ and
moreover, let $S'=S\setminus \{1\}$, then
$I^{1}_{init_{|S|-1}}(S)=1*I^{1}_{init_{|S'|-1}}(S')$.  Thus:
$$D_{K\sqcup L}(S)=D_{\Delta(Lk)}(S')=D_{lk(1,K)}(S')+D_{lk(1,L)}(S')$$
$$=D_{K}(S)+D_{L}(S)=D_{\Delta(K\dot{\cup}L)}(S).$$
The second and last equations are by Theorem \ref{prop.6} and
Proposition \ref{prop.5}. This completes the proof.$\square$

\subsection{How to shift a union over a simplex?}\label{KsimplexL}
In the case where $K\cap L=<\sigma>$ is a simplex (and all of its
subsets), we also get a formula for $\Delta(K \cup L)$ in terms of
$\Delta(K)$, $\Delta(L)$ and $\Delta(K\cap L)$. This case
corresponds to the topological operation called connected sum.

$Proof\  of\  Theorem\  \ref{clique-sum}$:  For a simplicial
complex $H$, let $\bar{H}$ denote the complete simplicial complex
$2^{H_{0}}$. The inclusions $H\hookrightarrow \bar{H}$ for
$H=K,L,<\sigma>$ induce a morphism from the commutative diagram
(\ref{diagram}) of $K$ and $L$ to the analogous commutative
diagram $(\bar{\ref{diagram}})$ of $\bar{K}$ and $\bar{L}$. By
functoriality of the sequence of the snake lemma, we obtain the
following commutative diagram:
\begin{equation} \label{snake*2}
\begin{CD}
0@> >>ker f@> >>ker g@> >>ker h@>\delta>>coker f@> >>... \\
@VV V @VVidV @VV V @VV V @VVidV \\
0@> >>ker \bar{f}@> >>ker \bar{g}@> >>ker \bar{h}@>\bar{\delta}>>coker \bar{f}@> >>... \\
\end{CD}
\end{equation}
where the bars indicate that ($\bar{\ref{diagram}}$) is obtained
from (\ref{diagram}) by putting bars over all the complexes and
renaming the maps by adding a bar over each map. Thus, if
$\bar{\delta}=0$ then also $\delta=0$, which, as we have seen,
implies (\ref{additive-formula}). The fact that
$\Delta(<\sigma>)=<\sigma>$ completes the proof.

We show now that $\bar{\delta}=0$. To simplify notation, assume
that $K$ and $L$ are complete complexes whose intersection is
$\sigma$ (which is a complete complex). Consider (\ref{diagram})
with $j=1$. (It is enough to prove Theorem \ref{clique-sum} for
$i=1$ as for every $i>1$, $S\subseteq [n]$ and $H$ a simplicial
complex on $[n]$, $I_{S}^{i}\cap H = \biguplus_{T\in
I_{S}^{i-1}}(I_{T}^{1}\cap H)$.) Let $m=m_{K}+m_{L}\in ker h$
where $supp(m_{K})\in K\setminus \sigma$. By commutativity of the
middle right square of (\ref{diagram}),
$\oplus_{A<_{L}S}f_{A}\lfloor(m_{K}),\oplus_{A<_{L}S}f_{A}\lfloor(m_{L})\in
\oplus_{A<_{L}S}\bigwedge^{1}(\sigma)$. If we show that
\begin{equation} \label{capIm}
\oplus_{A<_{L}S}f_{A}\lfloor(\bigwedge^{1+|S|}\sigma)=\oplus_{A<_{L}S}f_{A}\lfloor(\bigwedge^{1+|S|}K)
\bigcap \oplus_{A<_{L}S}\bigwedge^{1}(\sigma),
\end{equation}
then there exists $m'\in \bigwedge^{1+|S|}(\sigma)$ such that
$\oplus_{A<_{L}S}f_{A}\lfloor(m_{K}) =
\oplus_{A<_{L}S}f_{A}\lfloor(m') = f(m')$, hence
$\delta(m)=[f(m')]=0$ (where $[c]$ denotes the image of $c$ under
the projection onto $coker f$) i.e. $\delta=0$. (\ref{capIm})
follows from the intrinsic characterization of the image of the
maps it involves, given in Proposition \ref{intrinsicIm}. By
Proposition \ref{intrinsicIm}, the right hand side of
(\ref{capIm}) consists of all $x\in
\oplus_{R<_{L}S}\bigwedge^{1}K$ that satisfy $(a)$ and $(b)$ of
Proposition \ref{intrinsicIm} which are actually in
$\oplus_{A<_{L}S}\bigwedge^{1}(\sigma)$. By Proposition
\ref{intrinsicIm}, this is exactly the left hand side of
(\ref{capIm}). $\square$

The following generalizes a result of Kalai for graphs (\cite{56},
Lemma 3.7).
\begin{prop}\label{intrinsicIm}
Let $H$ be a complete simplicial complex with $H_{0}\subseteq
[n]$,  and let $S \subseteq [n]$, $|S|=s$. Then
$\oplus_{R<_{L}S}f_{R}\lfloor(\bigwedge^{1+s}H)$ is the set of all
$x=(x_{R}: R<_{L}S)\in \oplus_{R<_{L}S}\bigwedge^{1}H$ satisfying
the following:

(a) For all pairs $(i,R)$ such that $i\in R<_{L}S$:
 $<f_{i},x_{R}>=0$.

(b) For all pairs $(A,B)$ such that $A<_{L}S ,B<_{L}S$ and
$|A\vartriangle B|=2$: Denote $\{a\}=A\setminus B$ and
$\{b\}=B\setminus A$. Then $$-<f_{b},x_{A}>=(-1)^{sgn_{A\cup
B}(a,b)}<f_{a},x_{B}>$$ where $sgn_{A\cup B}(a,b)$ is the number
modulo $2$ of elements between $a$ and $b$ in the ordered set
$A\cup B$.
\end{prop}

$Proof$:
%Recall that
%$f_{R}=\sum_{S'\subseteq[n],|S'|=s}A_{RS'}e_{S'}$ where $A_{RS'}$
%is the minor of  the (generic) transition matrix $A=(\alpha_{ij})$
%w.r.t. the rows $R$ and columns $S'$. Let $T\subseteq H_{0}$,
%$|T|=s+1$. By (\ref{floorE}), $f_{R}\lfloor e_{T} = \sum_{t\in T}
%(-1)^{sgn(t,T)}A_{R,T\setminus \{t\}}e_{t}$
Let us verify first that every element in
$Im=\oplus_{R<_{L}S}f_{R}\lfloor(\bigwedge^{1+s}H)$ satisfies
$(a)$ and $(b)$. Let $y\in \bigwedge^{1+s}H$. If $i\in R<_{L}S$
then $<f_i,f_{R}\lfloor y>=<f_i\wedge f_{R},y>=<0,y>=0$, hence
$(a)$ holds. For $i\in T\subseteq [n]$ for some $n$, let
$sgn(i,T)=|\{t\in T: t<i\}|(mod\ 2)$. If $A,B<_{L}S$,
$\{a\}=A\setminus B$ and $\{b\}=B\setminus A$ then
$-<f_b,f_{A}\lfloor y>=-<f_b\wedge f_{A},y>=-(-1)^{sgn(b,A\cup
B)}<f_{A\cup B},y>=-(-1)^{sgn(b,A\cup B)}(-1)^{sgn(a,A\cup
B)}<f_a\wedge f_{B},y>=(-1)^{sgn_{A\cup B}(a,b)}<f_a,f_{B}\lfloor
y>$, hence $(b)$ holds.
% For every $i\in R$, by expanding in the
%$e_{j}$'s basis, we see that $<f_{i},f_{R}\lfloor e_{T}> =
%\sum_{t\in T} (-1)^{sgn(t,T)}\alpha_{it}A_{R,T\setminus \{t\}}$.
%Let $a(l,T)=(\alpha_{lt})_{t\in T}$ and let $A(l,R;T)$ be the
%matrix with first row $a(l,T)$ and below it the rows $a(r,T)$,
%$r\in R$, ordered by the order of their indices. Then
%$<f_{i},f_{R}\lfloor e_{T}> = det(A(i,R;T))=0$. Denote for short
%$Im=\oplus_{R<_{L}S}f_{R}\lfloor(\bigwedge^{1+s}H)$. Then every
%element of $Im$ satisfies $(a)$. Similarly, for $A,B<_{L}S$ such that
%$|A\vartriangle B|=2$ we compute $-<f_{b},f_{A}\lfloor
%e_{T}>=-det(A(b,A;T)=(-1)^{sgn_{A\cup
%B}(a,b)}det(A(a,B;T)=(-1)^{sgn_{A\cup B}(a,b)}<f_{a},f_{B}\lfloor
%e_{T}>$.

We showed that every element of $Im$ satisfies $(a)$ and $(b)$.
Denote by $X$ the space of all $x\in
\oplus_{R<_{L}S}\bigwedge^{1}H$ satisfying $(a)$ and $(b)$. It
remains to show that $dim(X)=dim(Im)$.

Following the proof of Proposition \ref{prop.*},
$dim(Im)=dim(\bigwedge^{1+s}H)-dim(\bigcap_{R<_{L}S}Ker_{s}f_{R}\lfloor(\bigwedge^{1+s}H)=
|\{T\in \Delta(H): |T|=s+1\}|-|\{T\in \Delta(H): |T|=s+1,
S_{(1)}^{(m)}\leq_{L}T\}|= |\{T\in \Delta(H): |T|=s+1,
init_s(T)<_{L}S \}|$. Let $h=|H_0|$ and $sum(T)=|\{ t\in T:
T\setminus \{t\}<_{L}S \}|$. Note that $\Delta(H)=2^{[h]}$.
Counting according to the initial $s$-sets, we conclude that in
case $s<h$,
\begin{eqnarray} \label{dimIm}
dim(Im)=|\{R: R<_{L}S, R\subseteq [h]\}|(h-s)- \nonumber\\
\sum\{sum(T)-1 : T\subseteq [h], |T|=s+1, init_{s}(T)<_{L}S,
sum(T)>1\}.
\end{eqnarray}
In case $s\geq h$, $dim(Im)=0$.

Now we calculate $dim(X)$. Let us observe that every $x\in X$ is
uniquely determined by its coordinates $x_R$ such that $R\subseteq
[h]$: Let $i\in R\setminus [h]$, $R<_{L}S$. Every $j\in
[h]\setminus R$ gives rise to an equation $(b)$ for the pair
$(R\cup j\setminus i,R)$ and every $j\in R\cap [h]$ gives rise to
an equation $(a)$ for the pair $(j,R)$. Recall that $x_R$ is a
linear combination of the form $x_R= \sum_{l\in
[h]}\gamma_{l,R}e_l$ with scalars $\gamma_{l,R}$. Thus, we have a
system of $h$ equations on the $h$ variables $(\gamma_{l,R})_{l\in
[h]}$ of $x_{R}$, with coefficients depending only on $x_{F}$'s
with $F<_{L}R$ (actually also $|F\cap [h]|=1+|R\cap [h]|$) and on
the generic $f_k$'s, $k\in [n]$. This system has a unique solution
as the $f_{k}$'s are generic. By repeating this argument we
conclude that $x_{R}$ is determined by the coordinates $x_F$ such
that $F\subseteq [h]$.

Let $x(h)$ be the restriction of $x\in X$ to its $\{x_R:
R\subseteq [h], R<_{L}S \}$ coordinates, and let $X(h)=\{x(h):
x\in X\}$. Then $dim(X(h))=dim(X)$.

Let $[a]$, $[b]$ be the matrices corresponding to the equation
systems $(a)$, $(b)$ with variables $(\gamma_{l,T})_{l\in
[h],T<_{L}S}$ restricted to the cases $T\subseteq [h]$ and
$A,B\subseteq [h]$, respectively. $[a]$ is an $s\cdot
|\{R\subseteq [h]: R<_{L}S\}|\times h\cdot |\{R\subseteq [h]:
R<_{L}S\}|$ matrix and $[b]$ is a $|\{A,B\subseteq [h]: A,B<_{L}S,
|A\vartriangle B|=2\}|\times h\cdot |\{R\subseteq [h]: R<_{L}S\}|$
matrix.

We observe that the row spaces of $[a]$ and $[b]$ have a zero
intersection. Indeed, for a fixed $R\subseteq [h]$, the row space
of the restriction of $[a]$ to the $h$ columns of $R$ is
$span\{f_{i}^{0}: i\in R\}$ (recall that $f_{i}^{0}$ is the
obvious projection of $f_i$ on the coordinates $\{e_j:j\in
H_0\}$), and the row space of the restriction of $[b]$ to the $h$
columns of $R$ is $span\{f_{j}^{0}: j\in [h]\setminus R\}$. But as
the $f_{k}^{0}$'s, $k\in [h]$, are generic, $span\{f_{k}^{0}: k\in
[h]\}=\bigwedge^{1}H$. Hence $span\{f_{i}^{0}: i\in R\}\cap
span\{f_{j}^{0}: j\in [h]\setminus R\}=\{0\}$. We conclude that
the row spaces of $[a]$ and $[b]$ have a zero intersection.
%if
%$v$ belongs to the row space of $[a]$ and $v_{R}\neq 0$ (recall
%that $v_R$ is a vector with $h$ coordinates) then $v_{R}\neq
%u_{R}$ for every $u$ in the row space of $[b]$, as the $f_{i}$'s
%are generic.

$[a]$ is a diagonal block matrix whose blocks are generic of size
$s\times h$, hence
\begin{eqnarray} \label{dim[a]}
rank([a])=s\cdot |\{R: R<_{L}S, R\subseteq [h]\}|
\end{eqnarray}
in case $s<h$.

Now we compute $rank([b])$. For $T\subseteq [h]$, $|T|=s+1$, let
us consider the pairs in $(b)$ whose union is $T$. If $(A,B)$ and
$(C,D)$ are such pairs, and $A\neq C$, then $(A,C)$ is also such a
pair. In addition, if $A,B,C$ are different (the union of each two
of them is $T$) then the three rows in $[b]$ indexed by $(A,B)$,
$(A,C)$ and $(B,C)$ are dependent; the difference between the
first two equals the third. Thus, the row space of all pairs
$(A,B)$ with $A\cup B=T$ is spanned by the rows indexed
$(init_{s}(T),B)$ where $init_{s}(T)\cup B=T$.

We verify now that the rows $\bigcup\{(init_{s}(T),B):
init_{s}(T)\cup B=T\subseteq [h], |T|=s+1, |B|=s\}$ of $[b]$ are
independent. Suppose that we have a nontrivial linear dependence
among these rows. Let $B'$ be the lexicographically maximal
element in the set of all $B$'s appearing in the rows $(A,B)$ with
nonzero coefficient in that dependence. There are at most $h-s$
rows with nonzero coefficient whose restriction to their $h$
columns of $B'$ is nonzero (they correspond to $A$'s with
$A=init_s(B\cup\{i\}$ where $i\in [h]\setminus B$). Again, as the
$f_{i}$'s are generic, this means that the restriction of the
linear dependence to the $h$ columns of $B'$ is nonzero, a
contradiction. Thus,
\begin{eqnarray} \label{dim[b]}
rank([b])=|\bigcup_{B<_{L}S}\{(init_s(T),B): init_{s}(T)\cup
B=T\subseteq [h], |T|=s+1, |B|=s\}|\nonumber\\
= \sum\{sum(T)-1 : T\subseteq [h], |T|=s+1, init_{s}(T)<_{L}S,
sum(T)>1\}.
\end{eqnarray}
(Note that indeed $B<_{L}S$ implies $init_{s}(T)<_{L}S$ as
$B\subseteq T$.)

For $s<h$, $dim(X(h))=h\cdot |\{R: R<_{L}S, R\subseteq
[h]\}|-rank([a])-rank([b])$, which by (\ref{dimIm}),
(\ref{dim[a]}) and (\ref{dim[b]}) equals $dim(Im)$. For $s\geq h$,
$dim(X(h))=0=dim(Im)$. This completes the proof. $\square$

As a corollary of Theorem \ref{clique-sum}, we get the following
combinatorial formula for shifting the union over a simplex of
simplicial complexes:
\begin{thm} \label{corKcupL}
%(our main result in this section).
 Let $K$ and $L$ be simplicial complexes where $K\cap
L=<\sigma>$ is a complete simplicial complex.
  Let $(K\cup L)_{0}=[n], [n]\supseteq T=\{ t_{1}<\dots<t_{j}<t_{j+1}\}$. Then: $$T\in\Delta(K\cup L)
  \Leftrightarrow t_{j+1}-t_{j}\leq D_{K}(T)+D_{L}(T)-D_{<\sigma>}(T).$$
In particular, any gluing of $K$ and $L$ along a $d$-simplex
results with the same shifted complex $\Delta(K\cup L)$, depending
only on $\Delta(K)$, $\Delta(L)$ and $d$.
\end{thm}
$Proof$: Put $i=1$ and $S=init_{|T|-1}(T)$ in Theorem
\ref{clique-sum}, and by Proposition \ref{prop.3} we are done.
 $\square$

\textbf{Remark}: For symmetric shifting, the analogous assertions
to Theorem \ref{clique-sum} and Theorem \ref{corKcupL} remain
true. For their proof one uses a symmetric variant of Proposition
\ref{intrinsicIm} (where condition $(a)$ is omitted, and condition
$(b)$ has a symmetric analogue).

\section{Shifting Near Cones}\label{sec4}
  A simplicial complex $K$ is called
a near cone with respect to a vertex $v$ if for every $j\in S\in
K$ also $S\cup v\setminus j\in K$. We are about to prove a
decomposition theorem for the shifted complex of a near cone, from
which the formula for shifting a cone (mentioned in the
introduction) will follow. As a preparatory step we introduce the
Sarkaria map, modified for homology.

\subsection{The Sarkaria map}\label{sec33}
    Let $K$ be a near cone with respect to a vertex $v=1$. Let
  $e=\sum_{i \in K_{0}}e_{i}$ and let $f=\sum_{i \in K_{0}}\alpha_{i}e_{i}$ be a linear combination
  of the $e_{i}$'s such that $\alpha_{i}\neq 0$ for every $i \in K_{0}$.
  Imitating the Sarkaria maps for cohomology \cite{sarkaria}, we get for homology
  the following linear maps:
  $$(\bigwedge K, e_{v}\lfloor) \overset{U}{\longrightarrow}(\bigwedge K, e\lfloor)
   \overset{D}{\longrightarrow}(\bigwedge K, f\lfloor)$$
 defined as follows: for $S\in K$
 $$U(e_{S})= \left\{ \begin {array}{ll}
   e_{S}-\sum_{i \in S}(-1)^{sgn(i,S)} e_{v\cup S\setminus i} & \textrm{if
   $v\notin S$}\\
   e_{S} & \textrm{if $v\in S$} \end {array} \right. $$

$$D^{-1}(e_{S})=(\prod_{i\in S}\alpha_{i})e_{S}.$$
It is justified to write $D^{-1}$ as all the $\alpha_{i}$'s are
non zero.
\begin{prop} \label{prop.9Sarc}
The maps $U$ and $D$ are isomorphisms of chain complexes. In
addition they satisfy the following 'grading preserving' property:
if $S\cup T \in K$, $S\cap T=\emptyset$, then $$U(e_{S}\wedge
e_{T})=U(e_{S})\wedge U(e_{T})\  and\ D(e_{S}\wedge
e_{T})=D(e_{S})\wedge D(e_{T}).$$

\end{prop}
%\begin{proof}
$Proof$: The check is straight forward.  First we check that $U$ and
$D$ are chain maps. Denote $\alpha_{S}=\prod_{i\in S}\alpha_{i}$.
For every $e_{S}$, where $S \in K$, $D$ satisfies
$$D\circ e\lfloor (e_{S})= D(\sum_{j\in S}(-1)^{sgn(j,S)}e_{S\setminus j})
= \sum_{j\in
S}(-1)^{sgn(j,S)}\frac{\alpha_{j}}{\alpha_{S}}e_{S\setminus j}$$
and
$$f\lfloor\circ D(e_{S})=f\lfloor(\frac{1}{\alpha_{S}}e_{S}) = \sum_{j\in
S}(-1)^{sgn(j,S)}\frac{\alpha_{j}}{\alpha_{S}}e_{S\setminus j}.$$
For $U$: if $v\in S$ we have
$$U\circ e_{v}\lfloor (e_{S})=U(e_{S\setminus v}) =e_{S\setminus
v}-\sum_{i\in S\setminus v}(-1)^{sgn(i,S\setminus v)}e_{S\setminus
i} = \sum_{j\in S}(-1)^{sgn(j,S)}e_{S\setminus j}.$$ The last
equation holds because $v=1$. Further,
$$e\lfloor\circ U(e_{S})= e\lfloor(e_{S})=\sum_{j\in S}(-1)^{sgn(j,S)}e_{S\setminus
j}.$$ If $v\notin S$ we have
$$U\circ e_{v}\lfloor (e_{S})=U(0)=0$$
and
$$e\lfloor\circ U(e_{S})= e\lfloor(e_{S})-e\lfloor(\sum_{j\in S}(-1)^{sgn(j,S)}e_{S\cup v\setminus
j})=$$ $$\sum_{j\in S}(-1)^{sgn(j,S)}e_{S\setminus j}- \sum_{i\in
S}(-1)^{sgn(i,S)}\sum_{t\in S\cup v\setminus i}(-1)^{sgn(t,S\cup
v\setminus i)}e_{S\cup v\setminus\{i,t\}}=$$
$$\sum_{j\in S}(-1)^{sgn(j,S)}e_{S\setminus j}(1-(-1)^{sgn(v,S\cup v\setminus
j)}) - \sum_{j,i\in S,i\neq j}(-1)^{sgn(i,S)}(-1)^{sgn(j,S\cup
v\setminus i)}e_{S\cup v\setminus \{i,j\}}.$$ In the last line,
the left sum is zero as $v=1$, and for the same reason the right
sum can be written as:
$$=\sum_{j,i\in S, i<j}((-1)^{sgn(i,S)+sgn(j,S\setminus i)}+(-1)^{sgn(j,S)+sgn(i,S\setminus
j)})e_{S\cup v\setminus \{i,j\}}.$$ As $i<j$, the $\{ij\}$
coefficient equals
$$(-1)^{sgn(i,S)+sgn(j,S)+1}+(-1)^{sgn(j,S)+sgn(i,S)}=0,$$
hence $e\lfloor\circ U(e_{S})= U\circ e_{v}\lfloor (e_{S})$ for every
$S \in K$. By linearity of $U$ and $D$ (and of the boundary maps),
  we have that $U,D$ are chain maps. To show that $U,D$ are onto,
  it is enough to show that each $e_{S}$, where $S \in K$, is in
  their image. This is obvious for $D$. For $U$: if $v\in S$ then
  $U(e_{S})=e_{S}$, otherwise $e_{S}=U(e_{S})+\sum_{i \in S}(-1)^{sgn(i,S)} e_{v\cup S\setminus i}$
  , which is a linear combination of elements in $Im(U)$, so $e_{S}\in
  Im(U)$ as well. Comparing dimensions, $U$ and $D$ are
  also 1-1.

We now show that $U$ 'preserves grading' in the described above
sense (for $D$ it is clear). For disjoint subsets of $[n]$ define
$sgn(S,T)=|\{(s,t)\in S\times T: t<s\}|(mod\ 2)$. Let $S,T$ be
disjoint sets such that $S\cup T\in K$.  By $S\cup T$ we mean the
ordered union of $S$ and $T$ (and similarly for other set unions).

case 1: $v\notin S\cup T$.
 $$U(e_{S})\wedge U(e_{T})=$$ $$
e_{S}\wedge e_{T}+\sum_{i\in S}(-1)^{sgn(i,S)}e_{S\cup v\setminus
i}\wedge e_{T}+\sum_{j\in T}(-1)^{sgn(j,T)}e_{S}\wedge e_{T\cup
v\setminus j}=$$
$$(-1)^{sgn(S,T)}(e_{S\cup T}
+\sum_{l\in S\cup T}(-1)^{sgn(l,S\cup T)}e_{S\cup T\cup v\setminus
l})=$$
$$U(e_{S}\wedge e_{T}),$$
where the middle equation uses the fact that $v=1$, which leads to
the following sign calculation:
$$(-1)^{sgn(i,S)}(-1)^{sgn(S\cup v\setminus i,T)}= (-1)^{sgn(i,S)+sgn(S\setminus
i,T)}=$$
$$(-1)^{sgn(i,S)+sgn(S,T)+sgn(i,T)}=(-1)^{sgn(S,T)}(-1)^{sgn(i,S\cup
T)}.$$

case 2: $v\in S\setminus T$.$$U(e_{S})\wedge U(e_{T})=e_{S}\wedge
(e_{T}-\sum_{t \in T}(-1)^{sgn(t,T)} e_{T\cup v\setminus t}
)=e_{S}\wedge e_{T}=U(e_{S}\wedge e_{T}).$$

case 3: $v\in T\setminus S$. A similar calculation to the one for
case 2 holds.
 $\square$
\\
\textbf{Remark}: The 'grading preserving' property of $U$ and $D$
extends to the case where $S\cap T\neq \emptyset$ ($S,T\in K$),
but we won't use it here. One has to check that in this case
(where clearly $e_{S}\wedge e_{T}=0$):
$$U(e_{S})\wedge U(e_{T})=D(e_{S})\wedge D(e_{T})=0.$$
%\end{proof}

\subsection{Shifting a near cone}
  %A simplicial complex $K$ is called
%a near cone with respect to a vertex $v$ if for every $j\in S\in
%K$ also $S\cup v\setminus j\in K$.

 We now prove a decomposition
theorem for the shifted complex of a near cone.

\begin{thm} \label{thm.10'near-cone}
 Let $K$ be a near cone on a vertex set $[n]$ with respect to a vertex $v=1$.
Let $X=\{f_{i}: 1\leq i\leq n\}$ be some basis of $\bigwedge^{1}
K$ such that $f_{1}$ has no zero coefficients as a linear
combination of some given basis elements $e_{i}$'s of
$\bigwedge^{1} K$, and such that for
$g_{i-1}=f_{i}-<f_{i},e_{1}>e_{1}$, $Y=\{g_{i}: 1\leq i\leq n-1\}$
is a linearly independent set.
 Then:
$$\Delta_{X}(K)=(1*\Delta_{Y}(lk(v,K)))\cup B$$
where $B$ is the set $\{S\in \Delta_{X}(K): 1 \notin S\}$.
\end{thm}
%\begin{proof}
Note that $1*\Delta_{Y}(lk(v,K))$ is not a simplicial complex, but
merely a collection of faces of $\Delta_{X}(K)$ which contain $1$.
We claim that this collection equals the set of all faces of
$\Delta_{X}(K)$ which contain $1$.
\\*
 $Proof$: Clearly for every $l\geq 0$:
 $Ker_{l}e_{v}\lfloor =
\bigwedge^{l+1}ast(v,K)$ and $Im_{l}e_{v}\lfloor =
\bigwedge^{l}lk(v,K)$. Using the Sarkaria map $D\circ U$ (see
Proposition \ref{prop.9Sarc}), we get that $Im f_{1}\lfloor$ is
isomorphic to $\bigwedge lk(v,K)$ and is contained (because of
'grading preserving') in a sub-exterior-algebra generated by the
elements
$b_{i}=DU(e_{i})=\frac{1}{\alpha_{i}}e_{i}-\frac{1}{\alpha_{v}}e_{v}$,
$i\in K_{0}\setminus v$.
% Let $f_{1}=f$ as
%above. Let $f_{i}$, $2\leq i \leq n=|K_{0}|$ be generic vectors
%with respect to $f$ and to all $e_{i}$'s
%% $b_{i}$'s
%so they are also generic linear combinations of the $b_{i}$'s.
 Let $S\subseteq [n], |S|=l, 1\notin
S$. Recall that $(g\wedge f)\lfloor h = g\lfloor(f\lfloor h)$. Now
we are prepared to shift.
$$\bigcap_{R<_{L}1\cup S}Ker_{l} f_{R}\lfloor \cong Ker_{l} f_{1}\lfloor \oplus
\bigcap_{1\notin R<_{L}S}Ker f_{R}\lfloor (Im_{l} f_{1}\lfloor
\rightarrow \mathbb{R}),$$ which by the Sarkaria map is isomorphic
to
%$$Ker_{l} f\lfloor \oplus
%\bigcap_{1\notin R<_{L}S}Ker f_{R}\lfloor (Im_{l} f\lfloor
%\rightarrow \mathbb{R})
\begin{equation}\label{eq3.th10}
\wedge^{l+1}ast(v,K)\oplus \bigcap_{1\notin T<_{L}S}Ker
((DU)^{-1}f_{T})\lfloor (\wedge^{l}lk(v,K)\rightarrow \mathbb{R}).
\end{equation}

Denote by $\pi_{t}$ the natural projection $\pi_{t}: span\{e_{R}:
|R|=t\}\rightarrow span\{e_{R}: |R|=t, v\notin R\}$. Then
$$Ker ((DU)^{-1}f_{T})\lfloor (\wedge^{l}lk(v,K)\rightarrow
\mathbb{R}) = Ker(\pi_{l}(DU)^{-1}f_{T})\lfloor
(\wedge^{l}lk(v,K)\rightarrow \mathbb{R}).$$ Now look at the
matrix $(<\pi_{l}\circ(DU)^{-1}f_{T} , e_{R}>)$ where $1\notin T
<_{L}S, R\in lk(v,K)_{l-1}$. It is obtained from the matrix
$(<g_{T-1} , e_{R}>)$, where $1\notin T <_{L}S, R\in
lk(v,K)_{l-1}$, by multiplying each column $R$ by the non-zero
scalar $\prod_{i\in R}\alpha_{i}$ (by $T+i$ we mean the set
$\{t+i: t\in T\}$). Thus, restricting to the first $m$ rows of
each of these two matrices we get matrices of equal rank, for
every $m$. This means, in terms of kernels and using the proof of
Proposition \ref{prop.2}, that in (\ref{eq3.th10}) we can replace
$(DU)^{-1}f_{T}$ by $g_{T-1}$ (note that the proof of Proposition
\ref{prop.2} can be applied to non-generic shifting as well). We
get (putting $Q=T-1$):
%The $l$-compound matrix $(B_{ST})$ of $B$ ($|S|=|T|=l$) is
%obtained from $(A_{ST})$ by multiplying column $T$ by the non zero
%scalar $\prod_{i\in T}\alpha_{i}$ for each $T\subseteq[n],|T|=l$.
%Restricting to the first $m$ rows of $(B_{ST})$ and of $(A_{ST})$
%we get matrices of equal rank, for every $m$. In terms of kernels,
%this means that in (\ref{eq3.th10}) we can replace
%$(DU)^{-1}f_{T}$ by $g_{T-1}$ (by $S+i$ we mean $\{s+i: s\in
%S\}$), and get:

$$dim \bigcap_{R<_{L}1\cup S}Ker_{l} f_{R}\lfloor(K) =
dim \wedge^{l+1}ast(v,K)+ dim \bigcap_{Q<_{L}S-1}Ker_{l-1}
g_{Q}\lfloor (lk(v,K)).$$
 As the left summand in the right hand side is a constant
independent of $S$, it is canceled when applying the last part of
Proposition \ref{prop.2}, and we get:
$$1\dot{\!\cup}S \in \Delta_{X}(K)\  \Leftrightarrow $$
$$dim \bigcap_{R<_{L}1\cup S}Ker_{l} f_{R}\lfloor(K)
> dim \bigcap_{R\leq_{L}1\cup S}Ker_{l} f_{R}\lfloor(K)\
\Leftrightarrow$$ $$dim \bigcap_{T<_{L}S-1}Ker_{l-1}\ g_{T}\lfloor
(lk(v,K)) > dim \bigcap_{T\leq_{L}S-1}Ker_{l-1}\
g_{T}\lfloor(lk(v,K))\ \Leftrightarrow$$
$$S-1\in\Delta_{Y}(lk(v,K)).$$

 Thus we get the
claimed decomposition of $\Delta_{X}(K)$. $\square$
\\
%\end{proof}
%Remark: we've already known that $\Delta(K)$ contains the cone in
%the above decomposition (because $L\subseteq K \Rightarrow
%\Delta(L)\subseteq \Delta(K)$), but that $B$ doesn't contain faces
%with the vertex '1' is new.

As a corollary we get the following decomposition theorem for the
generic shifted complex of a near cone.

\begin{thm} \label{thm.10near-cone}
 Let $K$ be a near cone with respect to a vertex $v$. Then
$$\Delta(K)=(1*\Delta(lk(v,K)))\cup B,$$
where $B$ is the set $\{S\in \Delta(K): 1 \notin S\}$.
\end{thm}
$Proof$: Apply Theorem\ref{thm.10'near-cone} for the case where $X$ is generic. In this
case, $Y$ is also generic, and the theorem follows .$\square$
%\end{proof}
%\subsection{Shifting cone}

As a corollary we get the following property \cite{skira}:
\begin{cor} \label{prop.11cone}
$\Delta\circ Cone = Cone \circ\Delta.$
\end{cor}
$Proof$: Consider a cone over $v$: $\{v\}*K$. By Theorem
\ref{thm.10near-cone}, $\{1\}*\Delta(K)\subseteq \Delta(\{v\}*K)$,
but those two simplicial complexes have equal $f$-vectors, and
hence, $\{1\}*\Delta(K)=\Delta(\{v\}*K)$.$\square$\\
%\end{proof}
\textbf{Remarks}: (1) Note that by associativity of the join
operation, we get by Corollary \ref{prop.11cone}:
$\Delta(K[m]*K)=K[m]*\Delta(K)$ for every $m$, where $K[m]$ is the
complete simplicial complex on $m$ vertices.

(2) Using the notation in Theorem \ref{thm.10'near-cone} we get:
$$\Delta_{X}\circ Cone = Cone \circ\Delta_{Y}.$$
%Remark:
%It is worth mentioning Duval's decomposition theorem for
%any simplicial complex in \cite{Du28}.
%(2) Looking closer at the proof of Theorem \ref{thm.10near-cone},
%we get the following generalization for more kinds of shifting:

\begin{de}\label{i near cone}
%Define:
 $K$ is an $i-near \ cone$ if there exist a sequence of simplicial
 complexes
$K=K(0)\supset K(1)\supset\dots \supset K(i)$ such that for every
$1\leq j\leq i$ there is a vertex $v_{j}\in K(j-1)$ such that
$K(j)=ast(v_{j},K(j-1))$ and $K(j-1)$ is a near cone w.r.t.
$v_{j}$.
\end{de}
\textbf{Remark}: An equivalent formulation is that there exists a
permutation $\pi: K_{0}=[n]\rightarrow [n]$ such that

$$\pi(i)\in S\in K, 1\leq
l<i \Rightarrow (S\cup \pi(l)\setminus \pi(i))\in K,$$ which is
more compact but less convenient for the proof of the
 following generalization of Theorem
\ref{thm.10near-cone}:
\begin{cor} \label{prop.12near-near-cone}
Let $K$ be an i-near cone. Then
$$\Delta(K)=B\cup\bigcup^{.}_{1\leq j\leq i}j*(\Delta(lk(v_{j},K(j-1)))+j),$$
where $B=\{S\in \Delta(K): S\cap[i]=\emptyset\}$.
\end{cor}
$Proof$:
%Let us prove this for $i=2$
The case $i=1$ is Theorem \ref{thm.10near-cone}.

By induction hypothesis,
$\Delta(K)=\tilde{B}\cup\dot{\!\bigcup}_{1\leq j\leq
i-1}j*(\Delta(K(j-1))+j)$ where $\tilde{B}=\{S\in \Delta(K):
S\cap[i-1]=\emptyset\}$. We have to show that
\begin{equation}\label{eq.th12}
\{S\in\Delta(K): min\{j\in
S\}=i\}=i*(\Delta(lk(v_{i},K(i-1)))+i)).
\end{equation}

For $|S|=l$ with $min\{j\in S\}=i$, we have
$$\bigcap_{R<_{L}S}Ker_{l-1} f_{R}\lfloor(K)\ =$$
\begin{equation}\label{eq.th12+}
(\bigcap_{j<i}\ \bigcap_{R:|R|=l, j\in R} Ker_{l-1}
f_{R}\lfloor)\cap (\bigcap_{R<_{L}S: min(R)=i}Ker_{l-1}
f_{R}\lfloor).\
\end{equation}

By repeated application of Proposition \ref{prop.1}, for each
$j<i$,  $$\bigcap_{R:|R|=l, j\in R} Ker_{l-1}
f_{R}\lfloor=Ker_{l-1}f_{j}.$$ Hence, (\ref{eq.th12+}) equals
$$(\bigcap_{j<i}Ker_{l-1} f_{j}\lfloor)\ \cap\ (\bigcap_{R<_{L}S: min(R)=i}Ker_{l-1}
f_{R}\lfloor)\ =$$
$$\bigcap_{R<_{L}S: min(R)=i}Ker_{l-1} f_{R}\lfloor(A),$$ where
$A=\bigcap_{j<i}Ker f_{j}\lfloor(\bigwedge K)$, taking the kernels in each
dimension.
% If we show that By
%Sarkaria map, $$\wedge(ast(v_{1},K))\cong Ker e_{v_{1}}\lfloor(K)
%\cong Ker f_{1}\lfloor(K)$$ as graded algebras. By assumption,
%$K(1)=ast(v_{1},K)$ is a near cone with respect to $v_{2}$, say.
%Applying Sarkaria map for $K(1)$, we get:
%$$(\wedge K(1), e_{v_{2}}\lfloor) \cong (\wedge K(1),f_{2}\lfloor)
%\cong (Ker f_{1}\lfloor, f_{2}\lfloor)$$

By repeated application of the Sarkaria map, we get that
$A\cong\bigwedge K(i-1)$ as 'graded' chain complexes. We will show
now that
\begin{equation}\label{eq.th12++}
dim \bigcap_{R<_{L}S: min(R)=i}Ker_{l-1} f_{R}\lfloor(A)= dim
\bigcap_{R<_{L}S-(i-1)}Ker_{l-1} f_{R}\lfloor(\bigwedge K(i-1)).
\end{equation}
Let $\varphi: \bigwedge K(i-1)\rightarrow A$ be the Sarkaria
isomorphism, and let $f$ be generic w.r.t. the basis
$\{e_{i},..,e_{n}\}$ of $\bigwedge^{0}K(i-1)$. Then $\varphi(f)$
is generic w.r.t. the basis $\{\varphi(e_{i}),..,\varphi(e_{n})\}$
of $A$. We can choose a generic $\bar{f}$ w.r.t.
$\{e_{1},..,e_{n}\}$ such that
$<\bar{f},\varphi(e_{j})>=<\varphi(f),\varphi(e_{j})>$ for every
$i\leq j\leq n$. Actually, we can do so for $n-i$ generic
$f_{j}$'s simultaneously (as multiplying a nonsingular matrix over
a field by a generic matrix over the same field results in a
generic matrix over that field). We get that
$$\bigcap_{R<_{L}S: min(R)=i}Ker_{l-1} \bar{f}_{R}\lfloor(A)= \bigcap_{R<_{L}S-(i-1)}Ker_{l-1} \varphi(f_{R})\lfloor(A)$$
$$\cong \bigcap_{R<_{L}S-(i-1)}Ker_{l-1} f_{R}\lfloor(\bigwedge K(i-1)).$$
As both the $f_{i}$'s and the $\bar{f}_{i}$'s are generic,
$\bigcap_{R<_{L}S: min(R)=i}Ker_{l-1} f_{R}\lfloor(A) \cong
\bigcap_{R<_{L}S: min(R)=i}Ker_{l-1} \bar{f}_{R}\lfloor(A)$ and
(\ref{eq.th12++}) follows. By applying Theorem
\ref{thm.10near-cone} to the near cone $K(i-1)$, we see that
(\ref{eq.th12}) is true, which completes the proof. $\square$
\\

From our last corollary we obtain a new proof of a well known
property of algebraic shifting, proved by Kalai \cite{D^2=D}:
\begin{cor} \label{prop.13 delta^2}
$\Delta^{2}=\Delta.$
\end{cor}
$Proof$: For every simplicial complex $K$ with $n$ vertices, $\Delta(K)$ is shifted,
(hence an $n$-near cone), and so are all the $lk(i,(\Delta K)(i-1))$'s associated to
it. By induction on the number of vertices, $\Delta(lk(i,(\Delta K)(i-1)))=lk(i,(\Delta
K)(i-1))-i$ for all $1\leq i\leq n$. Thus, applying Corollary
\ref{prop.12near-near-cone} to the $n$-near cone $\Delta(K)$, we get
$\Delta(\Delta(K))= \Delta(K)$.$\square$
%\end{proof}

 % What can be said about shifting of general join of simplicial
 % complexes will be discussed in the next section.

\section{Shifting Join of Simplicial Complexes}\label{sec k*l}
%\subsection{Shifting join}
Let $K,L$ be two disjoint simplicial complexes (they include the
empty set), and denote their join by $K*L$, i.e. $$K*L=\{S\cup T:
S\in K, T\in L\}.$$ $K*L$ is also a simplicial complex, and using
the K\"{u}nneth theorem with field coefficients (see \cite{Mun},
Theorem 58.8 and ex.3 on p.373) we can describe its homology in
terms of the homologies of $K$ and $L$:
$$H_{i}(K\times L)\cong \bigoplus_{k+l=i}H_{k}(K)\otimes
H_{l}(L)$$

and $$0\rightarrow \tilde{H}_{p+1}(K*L)\rightarrow
\tilde{H}_{p}(K\times L)\rightarrow
 \tilde{H}_{p}(K)\oplus \tilde{H}_{p}(L)\rightarrow 0.$$
Recalling that $\beta_{i}(K)=|\{S\in \Delta(K): |S|=i+1, S\cup
1\notin \Delta(K)\}|$ (Bj\"{o}rner and Kalai \cite{BK}), we get a
description of the number of faces in $\Delta(K*L)_{i}$ which
after union with $\{1\}$ are not in $\Delta(K*L)$, in terms of
numbers of faces of that type in $\Delta(K)$ and $\Delta(L)$. In
particular, if the dimensions of $K$ and $L$ are strictly greater
than $0$, the K\"{u}nneth theorem implies:
$$\beta_{dim(K*L)}(K*L)=\beta_{dim(K)}(K)\beta_{dim(L)}(L)$$
and hence
$$|\{S\in \Delta(K*L): 1\notin S, |S|=dim(K*L)+1\}|=$$
$$|\{S\in \Delta(K): 1\notin S, |S|=dim(K)+1\}|\times|\{S\in \Delta(L):
1\notin S, |S|=dim(L)+1\}|.$$ We now show that more can be said
about the faces of maximal size in $K*L$ that represents homology
of $K*L$, i.e. proving Theorem \ref{K*L,max}.
%\begin{thm}\label{K*L,max}
%Let $|(K*L)_{0}|=n$. For every $i\in [n]$ $$|\{S\in \Delta(K*L):
%[i]\cap S=\emptyset, |S|=dim(K*L)+1\}|=$$ $$|\{S\in \Delta(K):
%[i]\cap S=\emptyset, |S|=dim(K)+1\}| *$$ $$* |\{S\in \Delta(L):
%[i]\cap S=\emptyset, |S|=dim(L)+1\}|$$
%\end{thm}
\\

$Proof\  of\  Theorem\  \ref{K*L,max}$: For a generic
$f=\sum_{v\in K_{0}\cup L_{0}}\alpha_{v}e_{v}$ decompose
$f=f(K)+f(L)$ with supports in $K_{0}$ and $L_{0}$ respectively.
Denote $dim(K)=k, dim(L)=l$, so $dim(K*L)=k+l+1$. Observe that
$(f(K)\lfloor(K*L)_{k+l+1})\cap (f(L)\lfloor(K*L)_{k+l+1})=\{0\}$.
Denote by $f\lfloor(K)$ the corresponding generic boundary
operation on $span\{e_{S}:S\in K\}$, and similarly for $L$.
Looking at $\bigwedge(K*L)$ as a tensor product $(\bigwedge
K)\otimes (\bigwedge L)$ we see that $Ker_{k+l+1}f(K)\lfloor$
equals $Ker_{k}f(K)\lfloor\arrowvert_{\bigwedge K}\otimes
\bigwedge^{1+l} L$, and also $Ker_{k}f(K)\lfloor|_{\bigwedge
K}\cong Ker_{k}f\lfloor(K)$, and similarly when changing the roles
of $K$ and $L$. Hence, we get

$$Ker_{k+l+1}f\lfloor =$$ $$Ker_{k+l+1}f(K)\lfloor\ \cap
\ Ker_{k+l+1}f(L)\lfloor \cong$$ $$Ker_{k}f\lfloor(K)\otimes
Ker_{l}f\lfloor(L).$$
For the first $i$ generic $f_{j}$'s, by the
same argument, we have:
$$\bigcap_{j\in [i]}Ker_{k+l+1}f_{j}\lfloor =$$ $$\bigcap_{j\in [i]}Ker_{k+l+1}f_{j}(K)\lfloor\ \cap
\bigcap_{j\in [i]}Ker_{k+l+1}f_{j}(L)\lfloor \cong$$
$$\bigcap_{j\in [i]}Ker_{k}f_{j}\lfloor(K)\otimes \bigcap_{j\in [i]}Ker_{l}f_{j}\lfloor(L).$$
By Propositions \ref{prop.1} and \ref{prop.2} we get the claimed equation (\ref{join}).
$\square$
%end of proof
%Remark: Note that we did not use genericity (sp??) in Propositions
%\ref{prop.1}, \ref{prop.2} nor elsewhere in this proof, so for any
%transition matrix $A=(B|C)$ of $K*L$ where $B$,$C$ correspond to
%$K_{0},L_{0}$ respectively and are , we get:
%$$|\{S\in \Delta_{A}(K*L): [i]\cap S=\emptyset,
%|S|=dim(K*L)+1\}|=$$ $$|\{S\in \Delta_{B}(K): [i]\cap S=\emptyset,
%|S|=dim(K)+1\}| *$$ $$* |\{S\in \Delta_{R}(L): [i]\cap
%S=\emptyset, |S|=dim(L)+1\}|$$
\\
\textbf{Remark}: For symmetric shifting, the analogous assertion
to Theorem \ref{K*L,max} is false. As an example, let each of $K$
and $L$ consist of three points. Thus, $K*L=K_{3,3}$ is the
complete bipartite graph with $3$ vertices on each side. By
Theorem \ref{K*L,max}, $\{3,4\}\in \Delta(K_{3,3})$, but
$\{3,4\}\notin \Delta^{symm}(K_{3,3})$ where $\Delta^{symm}$
stands for the symmetric shifting operator (\cite{skira}, p. 128).

We now deal with the conjecture (\cite{skira}, Problem $12$)
\begin{equation}\label{eq**}
\Delta(K*L)=\Delta(\Delta(K)*\Delta(L)).
\end{equation}
We give a counterexample showing that it is false even if we
assume that one of the complexes $K$ or $L$ is shifted. Denote by
$\Sigma K$ the suspension of $K$, i.e. the join of $K$ with the
(shifted) simplicial complex consisting of two points.
\\
\textbf{Example}: Let $B$ be the graph consisting of two disjoint
edges. In this case $\Delta (\Sigma(B)) \setminus \Delta (\Sigma
(\Delta((B))) = \{\{1,2,6\}\}$ and $\Delta (\Sigma (\Delta(B)))
\setminus \Delta (\Sigma (B)) = \{\{1,3,4\}\},$ so (surprisingly)
we even get that
\begin{equation}\label{neq}
\Delta (\Sigma (B))<_{L}\Delta (\Sigma (\Delta(B))),
\end{equation}
where the lexicographic partial order on simplicial complexes is
defined (as in \cite{skira}) by: $K\leq_{L} L$ iff for all $r>0$
the lexicographically first $r$-face in $K\triangle L$ (if exists)
belongs to $K$.

\begin{conj}\label{suspension}
For any simplicial complex $K$: $\Delta (\Sigma (K))\leq_{L}\Delta
(\Sigma (\Delta(K)))$.
\end{conj}
This manuscript was first put on the math arXiv about two years ago. Very recently Satoshi Murai announced a proof of  Conjecture \ref{suspension}, and more generally, that $\Delta(K*L)\leq_{L}\Delta(\Delta(K)*\Delta(L))$ for any two simplicial complexes $K$ and $L$. 

\begin{conj}\label{suspension-topo.}(Topological invariance.)
Let $K_1$ and $K_2$ be triangulations of the same topological
space. Then $\Delta (\Sigma (K_1))<_{L}\Delta (\Sigma
(\Delta(K_1)))$ iff $\Delta (\Sigma (K_2))<_{L}\Delta (\Sigma
(\Delta(K_2)))$.
\end{conj}
It would be interesting to find out when equation (\ref{eq**})
holds.
%, or at least to give non-trivial sufficient conditions on
%$K,L$ that ensure (\ref{eq**}) holds.
If both $K$ and $L$ are shifted, it trivially holds as
$\Delta^{2}=\Delta$. By the remark to Corollary \ref{prop.11cone}
it also holds if $K$, say, is a complete simplicial complex.

\section*{Acknowledgments}  I would like to express my profound thanks to my supervisor
prof. Gil Kalai, for numerous helpful discussions. 
I deeply thank Isabella Novik for valuable remarks on early
versions of this manuscript. 
Thanks go also to the referees; one of them deserves special thanks for a very careful reading and very
helpful suggestions.

\end{document}